\documentclass[12pt]{article}
\usepackage{amssymb,amsmath}

\begin{document}
\date{}
\title{Rigidity of noncompact complete manifolds with harmonic curvature}

\author{Seongtag Kim}

\maketitle
\begin{center}
$\;$\\
$^\dagger$Department of Mathematics Education,\\ Inha University, \\
Incheon, 402-751 Korea \\
 \textrm{e-mail:}  \texttt{stkim@inha.ac.kr}
\end{center}

\begin{abstract}
Let $(M,g)$ be a  noncompact complete $n$-manifold with harmonic
curvature and positive Sobolev constant.  Assume that $L_2$ norms
of Weyl curvature  and traceless Ricci curvature are finite. We
prove that $(M,g)$ is Einstein if $n \ge 5$ and  $L_{n/2}$
norms of Weyl curvature and traceless Ricci curvature are small
enough.
\end{abstract}


\newtheorem{theorem}{Theorem}
\newtheorem{lemma}[theorem]{Lemma}
\newtheorem{cor}[theorem]{Corollary}
\newtheorem{rmk}[theorem]{Remark}
\def\bea{\begin{eqnarray}}
\def\eea{\end{eqnarray}}
\newcommand{\be}{\begin{equation}}
\newcommand{\ee}{\end{equation}}
\newcommand{\beq}{\begin{eqnarray}}
\newcommand{\eeq}{\end{eqnarray}}
\newcommand{\beqn}{\begin{eqnarray*}}
\newcommand{\eeqn}{\end{eqnarray*}}
\newcommand{\pa}{\partial}
\newcommand{\pxi}{ {\pa \over \pa x_i}}
\newcommand{\pxj}{ {\pa \over \pa x_j}}
\newcommand{\pxk}{ {\pa \over \pa x_k}}
\newcommand{\pxl}{ {\pa \over \pa x_l}}
\newcommand\nm{\nonumber}
\newcommand\ep{\epsilon}
\newcommand\ew{\epsilon'}

\newcommand{\ai}[1]{\inf_{u\in C^{\infty}(M)}\frac{\int\limits_{#1}{|\nabla
u|^2+C_n S \, u^2 \; dV_g} }
 {(\int\limits_{#1}{u^{ 2n/( n-2) }\;  dV_g)^{(n-2)/n}  }  } }\par
\newcommand{\aI}[1]{\inf\limits_{u\in C_0^{\infty}(M-B_r)}  \frac {\int
\limits_{#1}{|\nabla u|^2+C_n S \, u^2 \; dV_g} }
 { (\int\limits_{#1}{u^{ 2n/(n-2)}\;  dV_g)^{(n-2)/n}  }    }}
\newcommand{\ab}[2]{  \frac{  \int\limits_{#1} |\nabla  {#2}|^2+C_n S
{#2}^2 dV_g} {     (\int\limits_{#1}{{#2}^{ 2n/(n-2)}\;
dV_g)^{(n-2)/n}} }}

\section{Introduction} 
One of the important problems in Riemannian geometry is to understand metrics near  the Einstein space or the constant curvature space. 
The metric with harmonic curvature is a natural candidate for this study.  
Since the curvature tensor is decomposed into orthogonal parts, we 
 compare this metric with the Einstein metric and 
the constant curvature space, when  Weyl curvature and
traceless Ricci curvature are small in $L_{n/2}$ sense.   By this comparison, we obtain some $L_{n/2}$-type rigidity of  noncompact complete metrics with harmonic curvature in Theorem \ref{th1}-\ref{th5}.
\vskip 0.3 true cm

A Riemannian $n$-manifold  $(M,g)$ is a manifold
with harmonic curvature if the divergence of curvature tensor
$Riem$ vanishes, i.e.,\bea 0\ &=&(\delta Riem)_{jkl} =\nabla^i
R_{ijkl}, \label{bi1} \eea which is equivalent to
 \bea
0= \nabla_k R_{jl}-\nabla_l R_{jk}. \label{bi2} \eea By the
Bianchi identity, scalar curvature  is constant for $n \ge 3$. An 
important example of metrics with harmonic curvature is the Einstein
space. The parallel Ricci curvature space is another example of the manifold with harmonic curvature, but 
the converse does not hold (see Derdzinski \cite{der}). In the space of  connections,  the Yang-Mills connection corresponds to a metric with harmonic 
curvature (see Bourguignon\cite{bou}). Important properties of the manifold with harmonic curvature 
were surveyed by Besse \cite{bes}.  Einstein manifolds and manifolds with harmonic
curvature share many important properties.   A natural question is when a manifold of harmonic 
curvature must be Einstein or flat. 
Hebey and Vaugon proved $L_{n/2}$-type rigidity for compact manifolds with harmonic curvature and 
positive scalar curvature \cite{hv}.  
In this paper, we provide some rigidity properties for noncompact complete
manifold with
harmonic curvature. For the  rigidity proof, we use an elliptic
estimation for the Laplacian of curvature tensors.

\section{Manifolds with harmonic curvature}

Let $(M,g)$ be a noncompact complete Riemannian manifold  of
dimension $n \ge 3$ with scalar curvature $R$. The Sobolev
constant $Q(M, g)$ is defined by
$$Q(M, g)\equiv \inf_{0\neq u\in C^{\infty}_0 (M)}  { { {     \int_M |\nabla u|^2+{{(n-2)}\over {4(n-1)} } R_g u^2 dV_g}
                  \over {\left( \int_M |u|^{   {2n}/(n-2)} dV_g
                  \right)^{(n-2)/n}}} } . $$
 There are noncompact complete Riemannian manifolds of negative scalar curvature  with positive  Sobolev
constant.  For example, any simply connected locally conformally
flat manifold has positive  Sobolev constant  \cite{sy}.  Contrary to
noncompact case, the Sobolev constant  of a given compact manifold
is determined by the sign of scalar curvature (see \cite {au}). It
is also known that $Q(M, g) \ge 0$ if a noncompact complete
manifold $(M, g) $ has zero scalar curvature \cite{del}. In this
paper, we show that the rigidity of noncompact complete manifolds
with harmonic curvature  holds  not only  for nonnegative scalar
curvature metrics but also for negative scalar curvature metrics with
positive Sobolev constant and the dimension $n \ge 5$.
\vskip 0.3 true cm
Let $E_{ij}$ be the traceless Ricci tensor, i.e.,
 $E_{ij}= R_{ij}-{1 \over 4} R g_{ij}$, $|E|=|E_{ij}|$, $W$ the
 Weyl curvature.  In this paper, we use $c$ and $c'$ to
 denote some positive constant, which can be varied.

\begin{theorem}\label{th1} Let  $(M,g)$  be a
 noncompact complete  Riemannian $n$-manifold with harmonic curvature and $Q(M,
 g)>0$.
 Assume that $\int_M
   |W|^{2}+|E|^2 \,  dV_g$ is finite, and {\rm(A)} or {\rm (B)} holds:
   \begin{enumerate}
   \item [{\rm (A)}]
   Scalar curvature $R\ge 0$.
   \item [{\rm(B)}]
   Scalar curvature $R <0$ and $n\ge 5$.
   \end {enumerate}
Then there exists a small number $c_0$  such that if $\int_M
   |W|^{n/2}+|E|^{n/2} \,  dV_g
 \le c_0$, then $(M,g)$  is Einstein.

\end{theorem}

\noindent {\it Proof. \ } We need to prove that $|E_{ij}|=0$. To
simplify notations, we will work in an orthonormal frame.  
 The Laplacian of traceless Ricci tensor is;
\bea \triangle E_{ij}&=&\nabla_k \nabla_k E_{ij} \nm \\
&=&\nabla_k \nabla_i E_{kj}\label{ha11} \label{q14}\\
&=&\nabla_i \nabla_k E_{kj}+R_{kikm}E_{mj}+R_{kijm}E_{km} \label{e15}\\
&=& R_{im}E_{mj}+R_{kijm}E_{km}, \label{p5}\eea \noindent where (\ref{bi2})
is used in (\ref{q14}) and $\nabla_k E_{kj}=0$ is used in
(\ref{p5}).  Since

\bea R_{kijm}&=& W_{kijm}+{
R\over{n(n-1)}}(g_{kj}g_{im}-g_{km}g_{ij}) \nm \\
&{}& \quad +{1 \over{n-2}}   \left(
E_{kj}g_{im}+E_{im}g_{kj}-E_{km}g_{ij}-E_{ij}g_{km}\right), \eea

\bea E_{ij}\triangle E_{ij} &=& W_{kijm}E_{km} E_{ij}+{
R\over{n(n-1)}}g_{kj}g_{im}E_{km}  E_{ij} \nm\\
&{}& \quad +{1 \over{n-2}}\left(
E_{kj}g_{im}+E_{im}g_{kj}\right)E_{km}E_{ij}+R_{im}E_{mj}E_{ij} \\
&=&W_{kijm} E_{km} E_{ij}+{ R\over{n(n-1)}}|E|^2 +{n \over{n-2}}
{\rm tr } E^3 +{R \over n} {\rm tr} E^2 \\
&=& W_{kijm} E_{km} E_{ij}+{ R\over{(n-1)}}|E|^2 +{n \over{n-2}}
{\rm tr } E^3 \label{e110} .
\eea

\noindent Note that  \bea    W_{kijm} E_{km} E_{ij}
&\le& \sqrt{{n-2} \over{2(n-1)}}|W||E|^2 ,\label{q26} \\
 {\rm tr } E^3 &\le& {{n-2}\over \sqrt{n(n-1)} }|E|^3  \label{q56}\eea
 \noindent and \bea
 |\nabla |E||^2 &\le& {n\over{n+2}}|\nabla
 E|^2  \label{ke1} 
\eea
 because $E_{ij}$ is a traceless Codazzi tensor (see \cite{hv}).
\noindent  Inequalities  (\ref{q26}) and
(\ref{q56}) are proved by Huisken \cite{hu}.
From (\ref{e110}),
\bea |E|\triangle |E|&=& |\nabla E|^2-|\nabla|E||^2+ E_{ij} \triangle E_{ij} \nm \\
&\ge& { 2 \over n} |\nabla |E||^2 -\sqrt{{n-2}
\over{2(n-1)}}|W||E|^2-\sqrt{n\over{n-1}}|E|^3+{{R |E|^2}\over
{n-1}} , \label{e1}\eea
 where (\ref{q26}--\ref{ke1}) is used in (\ref{e1}).  Let $u=|E|$.  Multiplying a smooth compact supported
function $ \phi^2 $ to (\ref{e1}) and  integrating on $M$,  for any positive constant $\ep_1$, \bea &{}&
\int_M (-\ep_1 +1+ {2/ n})\phi^2 |\nabla u|^2-\ep_1^{-1}|u
\nabla
\phi |^2  dV_g \\&{}&\le  \int_M (1+ {2/n})\phi^2 |\nabla u|^2+ 2 \phi u \nabla \phi \cdot \nabla u \, dV_g \\
&{}&  \le \int_M \phi^2 \Big[\sqrt{{n-2}
\over{2(n-1)}}|W|u^2+\sqrt{n\over{n-1}}u^3-{R\over {n-1}}u^2 \Big]
\ dV_g . \label{e22q}\eea Using (\ref{e22q}), for any positive constant $\ep_2$,
\bea &{}& \int_M  |\nabla ( \phi u)|^2 + {{n-2}\over {4 (n-1)}} R u^2 \phi^2 dV_g\\
 &{}& \le \int_M  (1+ \ep_2)\phi ^2 |\nabla u |^2 + (1+ \ep_2^{-1} ) u^2 |\nabla \phi |^2 +{{n-2}\over {4 (n-1)}} R u^2 \phi^2dV_g\\
 &{}& \le \int_M   c_1 u^2 |\nabla \phi |^2 +c_2 R u^2 \phi^2 +{\rm A \ }
dV_g,  \label{e34}\eea

\noindent where \bea c_1&=&(1+\ep_2) (-\ep_1 +1+ {2 /n})^{-1}
\ep_1^{-1}+\ep_2^{-1}+1,
\\
c_2&=&\Big((n-2)/4 - (1+\ep_2) (-\ep_1 +1+ {2/n})^{-1}
\Big)/(n-1),
\\
c_3&=& (1+\ep_2) (-\ep_1 +1+ {2 /n})^{-1}, \\
 {\rm A}&= & c_3 \big(\sqrt {{n-2} \over{2(n-1)}}|W|u^2+\sqrt{n
\over{n-1}}u^3 \big) \phi^2. \eea  Note that $c_2$ is positive
when $n\ge 5$ and  $\ep_1$  and $\ep_2$ are sufficiently small.
For any $n$, $c_2$ is negative when $\ep_1$ is sufficiently close
to $1+{2/n}$. Therefore we can  choose $\ep_1$ and $\ep_2$ so that $c_2 $ is
positive if $n\ge 5$ and $R<0$, and $c_2 $ is negative if $R>0$
for any given $n$, which makes the second term of (\ref{e34}) 
non-positive. Using the Sobolev constant $\Lambda_0\equiv Q(M, g)$
and (\ref{e34}),

\bea \Lambda_0 \left(\int_M (\phi u)^{2n/(n-2)} dv_g
\right)^{n/(n-2)}
 &\le& \int_M  |\nabla ( \phi u)|^2 + {{n-2}\over {4 (n-1)}} R
u^2
\phi^2  dV_g \nm \\
& \le &\int_M   c_1 u^2 |\nabla \phi |^2  +{\rm A \ } dV_g.
\label{e388}\eea Note that
 \bea &{}&\int_M  {\rm A \ } dV_g \nm \\
 &\le&    c_3
   \sqrt {{n-2}
\over{2(n-1)}} \left( \int_M |W|^{n/2} dV_g\right)^{2/n} \left(\int_M (u\phi)^{2n/(n-2)} dV_g \right)^{n/(n-2)}  \\
& {}& \quad + c_3 \sqrt{n \over{n-1}} \left(\int_M u^{n/2} dV_g
\right)^{2/n} \left(\int_M (u\phi)^{2n/(n-2)} dV_g \right)^{(n-2)/n}.
\nm\eea Since $\int_M |E|^{n/2}+|W|^{n/2} dV_g$ is sufficiently
small, $ \int_M {\rm A \ } dV_g $ can be absorbed into left hand
side of (\ref{e388}). Therefore, there exists a constant $c'$
  such that \bea
c' \left(\int_M (\phi u) ^{2n/(n-2)} \,  dV_g \right)^{(n-2)/n} \le  \int_M
  |\nabla \phi|^2 u^2 \,  dV_g . \label{c41}
  \eea
\noindent Let $B_t=\{ x\in M| d(x, x_0) \le t\}$ for some fixed
$x_0 \in M$ and choose
$\phi$ as \bea \phi=\begin {cases}1 & \text{on} \  B_t,\\
0 & \text{on} \  M-B_{2t}, \label{e42}
\\
|\nabla \phi|\le {2 \over t} & \text{on} \  B_{2t}-B_t, \\
\end{cases} \eea and $0\le \phi \le 1$.
From (\ref{c41}) \bea c' \left(\int_M (u \phi )^{2n/(n-2)} \, dV_g
\right)^{n/(n-2)} &\le& {4\over{t^2}}\int_{B(2t)-B(t)} u^2  \,
dV_g  . \eea By taking $t\to \infty$, we have $u=0$ since $\int_M
|E|^2 dV_g $ is finite. Therefore $(M,g)$ is Einstein.

\vskip 0.5 true cm
 Since there is no noncompact complete Einstein metric
with positive scalar curvature, Theorem \ref{th1} implies that
there should be a lower bound for $\int_M
   |W|^{n/2}+|E|^{n/2} \,  dV_g $
if
  $\int_M |W|^2+ |E|^2 \, dV_g$ is finite.

\begin{theorem}\label{thob}
Let  $(M,g)$  be a
 noncompact complete  Riemannian $n$-manifold with harmonic curvature and positive Sobolev constant. 
 Assume that  scalar curvature $R$ is  positive and $\int_M
   |W|^{2}+|E|^2 \,  dV_g$ is finite.
Then there exists a  positive constant $c$  such that $\int_M
   |W|^{n/2}+|E|^{n/2} \,  dV_g
 \ge c$.
 \end{theorem}
 
\vskip 0.3 true cm 
Next we get a rigidity result without   $L_2$  finiteness of
 Weyl curvature  and traceless Ricci curvature. 
\begin{theorem}\label{th11} Let  $(M,g)$  be a
 noncompact complete  Riemannian $n$-manifold with harmonic curvature and $Q(M,
 g)>0$.
 Assume that $n\ge 4$ and   {\rm(A')} or {\rm (B')} holds:
   \begin{enumerate}
   \item [{\rm (A')}]
   Scalar curvature $R\ge 0$.
   \item [{\rm(B')}]
   Scalar curvature $R <0$ and $n\ge 6$.
   \end {enumerate}
Then there exists a small number $c_0$  such that if $\int_M
   |W|^{n/2}+|E|^{n/2} \,  dV_g
 \le c_0$, then $(M,g)$  is Einstein.

\end{theorem}
\noindent {\it Proof. \ } We need to prove that $|E_{ij}|=0$.  Let $u=|E|$. 
Multiplying a smooth compact supported
function $ \phi ^2 u^{-2+n/2}$ to (\ref{e1}) and  integrating on $M$,  for any positive constant $\ep_1$,
 \bea &{}&  16 \left( {n/ 2}+{2/ n}-1 -\ep_1\right) n^{-2} \int_M \phi^2 |\nabla u^{n/4}|^2 dV_g\\
 &{}& \   \le \int_M \ep_1^{-1}|u^{n/2}  \nabla\phi |^2  + \phi^2 \Big[\sqrt{{n-2}
\over{2(n-1)}}|W|u^{n/2} \nm \\ &{}& \quad +\sqrt{n\over{n-1}}u^{1+n/2} -{R\over {n-1}}u^{n/2} \phi^2]
\ dV_g \label{f22q}
\eea

Using (\ref{f22q}), for any positive constant $\ep_2$,
\bea &{}& \int_M  |\nabla ( \phi u^{n/4})|^2 + {{n-2}\over {4 (n-1)}} R u^{n/2} \phi^2 dV_g\\
 &{}& \le \int_M  (1+ \ep_2)\phi ^2 |\nabla u^{n/4} |^2 + (1+ \ep_2^{-1} ) u^{n/2} |\nabla \phi |^2 +{{n-2}\over {4 (n-1)}} R u^{n/2} \phi^2dV_g\\
 &{}& \le \int_M   c'_1 u^{n/2} |\nabla \phi |^2 +c'_2 R u^{n/2} \phi^2 +{\rm A' \ }
dV_g,  \label{f34}\eea

\noindent where \bea c'_1&=&(1+\ep_2) n^2 [ 16 ( n/2+2/n-1-\ep_1)]^{-1}
\ep_1^{-1}+\ep_2^{-1}+1,
\\
c'_2&=&\Big((n-2)/4 - (1+\ep_2) n^2 [ 16 ( n/2+2/n-1-\ep_1)]^{-1}      
\Big)/(n-1),
\\
c'_3&=&(1+\ep_2) n^2 [ 16 ( n/2+2/n-1-\ep_1)]^{-1}, \\
 {\rm A'}&= & c'_3 \big(\sqrt {{n-2} \over{2(n-1)}}|W|u^{n/2}+\sqrt{n
\over{n-1}}u^{1+n/2} \big) \phi^2. \eea
Note that $c'_2$ is positive
when $n\ge 6$ and  $\ep_1$  and $\ep_2$ are sufficiently small.
For any $n$, $c'_2$ is negative when $\ep_1$ is sufficiently close
to $-1+{2/n}+{n/2}$. Therefore we can  choose $\ep_1$ and $\ep_2$ so that $c'_2 $ is
positive if $n\ge 6$ and $R<0$, and $c'_2 $ is negative if $R>0$
for any given $n$, which makes the second term of (\ref{f34}) 
non-positive. 
Using the Sobolev constant and (\ref{f34}),

\bea \Lambda_0 \left(\int_M (\phi u^{n/4})^{2n/(n-2)} dv_g
\right)^{n/(n-2)}
 &\le& \int_M  |\nabla ( \phi u^{n/4})|^2 + {{n-2}\over {4 (n-1)}} R
u^{n/2}
\phi^2  dV_g \nm \\
& \le &\int_M   c'_1 u^{n/2} |\nabla \phi |^2  +{\rm A' \ } dV_g.
\label{f388}\eea 
Note that
 \bea &{}&\int_M  {\rm A' \ } dV_g \nm \\
 &\le&    c'_3
   \sqrt {{n-2}
\over{2(n-1)}} \left( \int_M |W|^{n/2} dV_g\right)^{2/n} \left(\int_M (u^{n/4}\phi)^{2n/(n-2)} dV_g \right)^{n/(n-2)}  \\
& {}& \quad + c'_3 \sqrt{n \over{n-1}} \left(\int_M u^{n/2} dV_g
\right)^{2/n} \left(\int_M (u^{n/4}\phi)^{2n/(n-2)} dV_g \right)^{(n-2)/n}.
\nm\eea Since $\int_M |E|^{n/2}+|W|^{n/2} dV_g$ is sufficiently
small, $ \int_M {\rm A \ } dV_g $ can be absorbed into left hand
side of (\ref{f388}). 
 Therefore, there exists a constant $c'$
  such that \bea
c' \left(\int_M (\phi u^{n/4}) ^{2n/(n-2)} \,  dV_g \right)^{(n-2)/n} \le  \int_M
  |\nabla \phi|^2 u^{n/2 }\,  dV_g . \label{f41}
  \eea
Choosing a  compact supported function similar to (\ref{e42}),
we can easily show that $u^{n/4}=0$ on $M$. We conclude that $(M, g)$ is
a Einstein space.

\section{Rigidity of Einstein spaces} 
In this section, we provide an $L_{n/2}$-type rigidity for noncompact complete Einstein manifolds. 
 Previously, this type of rigidity was proved for compact Einstein manifolds with positive 
scalar curvature  \cite{sh, si, hv, it}. 
Throughout this section $(M, g)$ denotes  a noncompact complete
Einstein manifold of dimension $n\ge 4$. 

\begin{theorem}\label{th2} Let  $(M,g)$  be a
 noncompact complete  Einstein $n$-manifold with $Q(M,
 g)>0$ and  finite $\int_M
   |W|^{2} \,  dV_g$. Assume that $(M, g)$ satisfies {\rm(C)} or {\rm (D)}:
   \begin{enumerate}  
   \item[{\rm(C)}] $(M,g)$ is Ricci-flat.
   \item[{\rm(D)}] Scalar curvature $R $  is negative and  $n\ge 8$.
   \end{enumerate}    
 Then there exists a small constant $c$  such that if $\int_M
   |W|^{n/2}\,  dV_g
 \le c$, then $(M,g)$  is a constant curvature space.

\end{theorem}
\noindent {\it Proof. \ } We need to prove that $|W|=0$.  
For an Einstein metric, the Laplacian of Weyl tensor  is \bea
\Delta W_{ijkl}&=& \nabla_t \nabla_t  W_{ijkl} \\&=& {2 \over n} R
W_{ijkl}+ W*W , \label{w1}\eea where $W*W$  denotes quadratic terms of Weyl tensor (see Singer\cite{si}).
 Multiplying $ W_{ijkl}$ on
(\ref{w1}), \bea W_{ijkl} \Delta W_{ijkl} =
 {2 \over n} R \  |W_{ijkl}|^2+ W*W*W. \label{w2} \eea
 In above, $W*W*W$  denotes cubic terms of
Weyl tensor, which is bounded by $c_4 |W|^3$ for a positive
constant $c_4$. Multiplying a smooth compact supported
function $ \phi $ to (\ref{w2}) and  integrating on $M$, 
\bea \int_M (1-\ep_1) \phi^2 |\nabla W|^2 dV_g \le \int_M -{ 2
\over n}R^2 |W|^2 \phi^2 +c_4  \phi^2 |W|^3 + { 1\over
{\ep_1}}|\nabla \phi |^2 |W|^2 dV_g, \label{q47} \eea for any positive constant $\ep_1$. The following refined Kato
inequality was proved for Weyl tensor of the Einstein manifold, \bea
|\nabla |W||^2 \le {{n-1}\over{n+1}} |\nabla W|^2, \label{kaw}
\eea when $|W|\neq 0$  (see Bando et al. \cite{bkn}). Let $f=|W|$. For any positive constant $\ep_2$,
\bea |\nabla ( \phi f)|^2 &\le& (1+{1 \over \ep_2})| f \nabla
\phi|^2 +(1+\ep_2) | \phi \nabla f|^2
\\
&\le& (1+{1 \over \ep_2})| f \nabla \phi|^2
+(1+\ep_2){{n-1}\over{n+1}} | \phi \nabla W|^2. \label{q50}\eea
From (\ref{q47}), \bea  \Lambda_0 \left(\int_M (\phi
f)^{2n/(n-2)} dv_g \right)^{{n \over{n-2}}}
 &\le& \int_M  |\nabla ( \phi f)|^2 + {{n-2}\over {4 (n-1)}} R
f^2 \phi^2  dV_g \label{s1}\\
&\le& \int_M c_5 Rf^2 \phi^2+c_6 |\nabla \phi|^2 f^2 +c_7 \phi^2
f^3 dV_g  \label{q38},\eea  \noindent where \bea c_5 &=
&-{{2(1+\ep_2)(n-1)}\over
{n (n+1)(1-\ep_1)}}+{ {n-2}\over{4(n-1)}}, \\
 c_6 &=&1+{\ep_2}^{-1}+{{(1+\ep_2)(n-1)}\over {(n+1)(1-\ep_1)
\ep_1}}, \\
c_7&=&{{(1+\ep_2)(n-1)}\over {(n+1)(1-\ep_1) }}c_4. \eea Note that
$c_5$ is positive when $n\ge 8$ and $\ep_1$  and  $\ep_2$ are
sufficiently small. $c_5$ is negative when $\ep_1$ is sufficiently
close to $1$ for any given $n$. Therefore we can choose $\ep_1$
and $\ep_2$ so that  $c_5$ is positive when $R<0$ and $n\ge 8$, which makes the first term of (\ref{q38}) 
non-positive. The last term of (\ref{q38}) is bounded by \bea c_7\int_M \phi^2 f^3
dV_g \le c_7 \left( \int_M f^{n/2} dV_g\right)^{2/n} \left(\int_M
(f\phi)^{2n/(n-2)} dV_g \right)^{n/(n-2)}. \label{eq6}\eea Since
$\int_M |W|^{n/2} dV_g$ is sufficiently small, (\ref{eq6}) can be
absorbed into left hand side of (\ref{s1}). Therefore, there
exists a constant $c'$
  such that \bea
c' \left(\int_M (\phi f)^{2n/(n-2)} \,  dV_g \right)^{(n-2)/n} \le  \int_M
  |\nabla \phi|^2 f^2 \,  dV_g . \label{c24}
  \eea
Choosing a  compact supported function similar to (\ref{e42}),
we can easily show that $f=0$ on $M$. We conclude that $(M, g)$ is
a constant curvature space.

\vskip 0.3 true cm
From Theorem \ref{th1} and Theorem \ref{th2},
\begin{theorem}\label{th5} Let  $(M,g)$  be a
 noncompact complete  Riemannian $n$-manifold with harmonic curvature and $Q(M,
 g)>0$.
 Assume that $\int_M
   |W|^{2}+|E|^2 \,  dV_g$ is finite, and {\rm(E)} or {\rm (F)} holds: 
   \begin{enumerate}
    \item[{\rm(E)}]   Scalar curvature $R=0$.
   \item[{\rm(F)}] Scalar curvature $R $  is negative and $n\ge 8$. 
   \end {enumerate}
 Then there exists a small number $c$  such that if $\int_M
   |W|^{n/2}+|E|^{n/2} \,  dV_g
 \le c$, then $(M,g)$  is a constant curvature space.

\end{theorem}

\noindent{\bf Remarks.}  \   The constants $c_0$ and $c$ in Theorem
\ref{th1}-\ref{th5}, depend on the Sobolev constant. For the  compact manifold, these constants depend only on 
the dimension of a given manifold, because of the finiteness of volume \cite{si, hv, it}. The applicability of Theorem 
\ref{th1}-\ref{th5} to lower dimensional manifolds needs further study.

\end{document}